\documentclass[conference]{IEEEtran}

\usepackage{amsmath,amssymb,amscd}
\usepackage{cite}
\usepackage{graphicx}

\def\e{\epsilon}

\def\bp{\begin{proposition}}
\def\ep{\end{proposition}}
\def\bt{\begin{theo}}
\def\et{\end{theo}}
\def\be{\begin{equation}}
\def\ee{\end{equation}}
\def\bl{\begin{lemma}}
\def\el{\end{lemma}}
\def\bc{\begin{corollary}}
\def\ec{\end{corollary}}
\def\pr{\noindent{\bf Proof: }}

\def\bd{\begin{definition}}
\def\ed{\end{definition}}
\def\C{{\mathbb C}}

\def\min{{\rm min\,}}
\def\max{{\rm max\,}}

\def\min{{\rm min\,}}
\def\max{{\rm max\,}}

\newtheorem{theo}{Theorem}[section]
\newtheorem{lemma}{Lemma}[section]
\newtheorem{definition}{Definition}[section]
\newtheorem{corollary}{Corollary}[section]
\newtheorem{proposition}{Proposition}[section]

\begin{document}
	
\title{Accuracy of spike-train Fourier reconstruction for colliding nodes}

\author{\IEEEauthorblockN{Andrey Akinshin\IEEEauthorrefmark{1},
Dmitry Batenkov\IEEEauthorrefmark{2},
Yosef Yomdin\IEEEauthorrefmark{3},
}
\IEEEauthorblockA{\IEEEauthorrefmark{1}\IEEEauthorrefmark{3}Department of Mathematics,\\
The Weizmann Institute of Science, Rehovot 76100, Israel}
\IEEEauthorblockA{\IEEEauthorrefmark{1}Laboratory of Inverse Problems of Mathematical Physics,\\
Sobolev Institute of Mathematics SB RAS, Novosibirsk 630090, Russia}
\IEEEauthorblockA{\IEEEauthorrefmark{2}Department of Computer Science,\\
Technion --- Israel Institute of Technology, Haifa 32000, Israel}
\IEEEauthorblockA{\IEEEauthorrefmark{1} Email: andrey.akinshin@weizmann.ac.il}
\IEEEauthorblockA{\IEEEauthorrefmark{2} Email: yosef.yomdin@weizmann.ac.il}
\IEEEauthorblockA{\IEEEauthorrefmark{3} Email: batenkov@cs.technion.ac.il}
}

\maketitle

\begin{abstract}
We consider signal reconstruction problem for signals $F$ of the form $ F(x)=\left(x\right)=\sum_{j=1}^{d}a_{j}\delta\left(x-x_{j}\right),$
from their Fourier transform ${\cal F}(F)(s)=\int_{-\infty}^\infty F(x)e^{-isx}dx.$ We assume ${\cal F}(F)(s)$ to be
known for each $s\in [-N,N],$ with an absolute error not exceeding $\e > 0$. We give an {\it absolute lower bound}
(which is valid with any reconstruction method) for the ``worst case'' error of reconstruction of $F$ from ${\cal F}(F),$
in situations where the nodes $x_j$ are known to form an $l$ elements cluster of a size $h \ll 1$. Using ``decimation''
algorithm of \cite{batenkov_numerical_2014,batenkov_accurate_2014} we provide an upper bound for the reconstruction error,
essentially of the same form as the lower one.  Roughly, our main result states that for $h$ of order $\frac{1}{N}\e^{\frac{1}{2l-1}}$
the worst case reconstruction error of the cluster nodes is of the same order $\frac{1}{N}\e^{\frac{1}{2l-1}}$, and hence the inside
configuration of the cluster nodes (in the worst case scenario) cannot be reconstructed at all. On the other hand, decimation
algorithm reconstructs $F$ with the accuracy of order $\frac{1}{N}\e^{\frac{1}{2l}}$.
\end{abstract}

\IEEEpeerreviewmaketitle


\section{Introduction}\label{Intro}
\setcounter{equation}{0}

In this paper, we provide lower and upper bounds for the Fourier reconstruction error, in the presence of noise, of spike-train signals in
the case of ``almost colliding'' (or clustering) nodes. The lower bound is obtained via the analysis of the behavior of the Fourier transform
under perturbation of the nodes and amplitudes in the cluster, and so it is valid (in the worst case scenario) for any reconstruction method.
The upper bound follows from the accuracy analysis of ``decimation'' reconstruction algorithm, as given in \cite{batenkov_numerical_2014,batenkov_accurate_2014}.

We hope that our analysis may clarify some aspects of the ``Super-resolution problem'' for spike-train signals with clustering nodes, as it
appears in many old and recent publications on the subject (see, as a very small sample, \cite{donoho_superresolution_1992,Don1,Lev.Ful,McC,Ode.Bar.Pis,Min.Kaw.Min}, recent publications
\cite{azais_spike,batenkov_numerical_2014,batenkov_accurate_2014,Bat.Sar.Yom,Bat.Yom1,Bat.Yom2,candes_towards_2014,
candes_super-resolution_2013,demanet_super-resolution_2013,fernandez-granda_support_2013,heckel_super-resolution_2014,liao_music_2014,moitra_threshold_2014,Dem.Ngu,Mor.Can,Bat.Yom.Sampta13,duval_exact_2013,Yom2},
and references therein).\par

Let us assume that the signal $F(x)$ is a spike-train, i.e. it is a priori known to be a linear combination of $d$ shifted $\delta$-functions:
\be \label{eq:equation_model_delta}
F(x)=F_{A,X}\left(x\right)=\sum_{j=1}^{d}a_{j}\delta\left(x-x_{j}\right),
\ee
where $A=(a_1,\ldots,a_d) \in {\mathbb R}^d, \ X=X_d=(x_1,\ldots,x_d) \in {\mathbb R}^d.$ We shall always assume that
$x_1\leq x_2\leq \ldots \leq x_d$. As for the measurements, we assume that the Fourier transform

\be\label{Four.Transf}
{\cal F}(F)(s)=\int_{-\infty}^\infty F(x)e^{-isx}dx
\ee
is known for each $s\in [-N,N],$ with an absolute error not exceeding $\e > 0$. So our input measurement is a function $\Phi(s)$ satisfying
$|\Phi(s)-{\cal F}(F)(s)|\leq \e$ for $s\in [-N,N]$.\par

The first goal of the present paper is to study the ``worst case'' accuracy of reconstruction of $F$ from $\Phi(s)$ in situations where the
nodes $x_j$ are known to form a cluster of a size $h \ll 1$, while being near-uniformly positioned inside the cluster. We give an {\it absolute
lower bound} for the reconstruction error of the nodes $x_j$ from the measured function $\Phi(s),$ which is valid independently of the
reconstruction method applied.\par

Our second goal is to give an upper bound for the reconstruction error of the nodes $x_j$, under the same assumptions as above. We show that
the decimation algorithm, combined with a homotopy continuation solving of the resulting algebraic equations, as described in \cite{batenkov_numerical_2014,batenkov_accurate_2014},
produces an error of essentially the same order of magnitude as the lower bound.\par

Shortly, our main result is as follows:\par

\noindent 1. {\it If certain $l$ nodes of $F$ form a cluster of a size $h \sim \frac{1}{N}\e^{\frac{1}{2l-1}}$, while being near-uniformly
positioned inside the cluster, then the worst case reconstruction error $\Delta$ of the cluster nodes is at least $Ch$.}\par

\noindent 2. {\it If for the same signal $F$ the measurements error is smaller than $\e_1\sim \e^{\frac{2l}{2l-1}}$ then the decimation algorithm
reconstructs the cluster nodes with the error $\Delta$ being at most $ch, \ c\ll C.$}\par

The ``practical'' conclusion could be that the inside configuration of the cluster nodes {\it cannot be reconstructed at all} from the Fourier
transform ${\cal F}(F)(s),$ known with the error $\e,$ for $s\in [-N,N],$ if the cluster size $h$ is smaller than  $\frac{1}{N}\e^{\frac{1}{2l-1}}$.
However, slightly reducing (to $\e_1$) the allowed magnitude of the measurements error, we can accurately and robustly reconstruct the cluster nodes
via the decimation algorithm.\par

The reconstruction error $\Delta \sim \frac{1}{N}\e^{\frac{1}{2l-1}}$ would make practical reconstruction of two colliding nodes very difficult,
and of three or more virtually impossible. However, our bound is the worst case one, and one can hope that for a random noise a typical
reconstruction accuracy may be much better.\par

Let us stress that our result is pretty close to the main result of \cite{donoho_superresolution_1992}, where Fourier sampling of atomic measures on non-uniform grids is
studied. In particular, the connection of the form $\e=C(Nh)^{2l-1}$ between the noise, the bandwidth, and the clustering geometry which can be
stably recovered, appears also in \cite{donoho_superresolution_1992}. Very recently similar bounds were obtained for superresolution of positive sources
in \cite{Mor.Can}, and for a Fourier recovery of sparse vectors in \cite{Dem.Ngu}. There are also apparent similarities with the classical result of Slepian in \cite{Sle}. 
Compare a discussion in \cite{candes_towards_2014} of the role of sparsity and clustering, as they appear in the superresolution problem, and, in particular, the discussion
in Sections 1.7 and 3.2 of \cite{candes_towards_2014} of the ``absolute lower bounds'' for the reconstruction error. We plan to further investigate the above connections.

\section{Main result}\label{Main}

To state our results we have to make some ``normalizing'' assumptions on the signal $F$ to be recovered. Indeed, if some amplitudes $a_j$ are small,
the reconstruction accuracy of the corresponding nodes drops, while larger $a_j$ imply higher accuracy. So we shall assume that the amplitudes
$A=(a_1,\ldots,a_d)$ of the signal $F$ satisfy the following assumption $A(m,M)$:

$$
0< m \leq |a_j|\leq M<\infty, \ j=1,\ldots,d.
$$

\bd\label{cluster} A signal $F_{A,X}$ as given by (\ref{eq:equation_model_delta}) is said to form an $(l,h,\rho)$-cluster $X$ if there is
an interval $I \subset {\mathbb R}$ of length $h$ which contains exactly $l$ nodes $X_l=\{x_\kappa, x_{\kappa+1}, \ldots, x_{\kappa+l-1}\}$
of $F$, while the minimal distance between the nodes in $X_l$ is at least $\rho h$, \ $\rho >0$.
\ed
\bd\label{cluster.dist}
For two ordered subsets $V=(v_1,\ldots,v_q)$, and $W=(w_1,\ldots,w_q)$ in ${\mathbb R}$ the distance $d(V,W)$ is defined as

\begin{multline*}
d(V,W)=\max_{s=1}^q |v_s-w_s|=||v_s-w_s||_{l^\infty} \geq \\
\geq \frac{1}{\sqrt q}||v_s-w_s||_{l^2}.
\end{multline*}
\ed

The following theorem is the first main result of the paper:

\bt\label{main}
Let a signal $F^0=F_{A^0,X^0},$ satisfying assumption $A(m,M)$, form an $(l,h,\rho)$-cluster
$X^0_l=\{x^0_\kappa, x^0_{\kappa+1}, \ldots, x^0_{\kappa+l-1}\}.$ Then there exist parameters $A^1, X^1$, satisfying assumption $A(\frac{m}{2},2M)$,
such that the distance $d(X^0_l,X^1_l)$ between $X^0_l$ and $X^1_l=\{x^1_\kappa, x^{1}_{\kappa+1}, \ldots, x^{1}_{\kappa+l-1}\}$ is at least $C_1h$,
while for $F^1=F_{A^1,X^1},$ and for each $s \in {\mathbb R}$ with $|s|\leq \frac{1}{2\pi h}$ we have

\be\label{Fourier.error}
|{\cal F}(F^0)(s)-{\cal F}(F^1)(s)|\leq C_2 (hs)^{2l-1}.
\ee
In particular, for $s\in [-N,N], \ N \leq \frac{1}{2\pi h},$ this difference does not exceed $C_2(hN)^{2l-1}.$ Here the constants $C_1$ and $C_2$
depend only on $m,M,l,\rho.$
\et
\par

The proof of Theorem \ref{main} is given in Section \ref{Proof} below. From this result we immediately deduce the following:

\bc\label{main.cor}
Assume that the noise $N(s)$ in the Fourier sampling ${\cal F}(s)$ on $[-N,N]$ may be an arbitrary function with the only restriction that
$|N(s)|\leq \e,$ where $0<\e\ll 1$, and put $h_\e=\frac{1}{N}(\frac{\e}{C_2})^{\frac{1}{2l-1}}.$ Let $F^0$ be any signal, satisfying assumption $A(m,M)$,
and forming an $(l,h_\e,\rho)$-cluster $X^0_l$. Let $F^1$ be the new signal produced from $F^0$ as in Theorem \ref{main}. Then for any reconstruction
algorithm $\cal R$ the worst-case error in reconstruction of the nodes of either $X^0_l$, or of $X^1_l$ is not smaller than $\frac{1}{2}C_1h_\e$.
\ec
\pr
We pick an ``adversary'' noise $N(s)$ to be identically zero for the sampling of $F^0$ and to be equal to $N(s)={\cal F}(F^0)(s)-{\cal F}(F^1)(s)$ for
the sampling of $F^1$ In both cases, by our choice of $h_\e$ and by Theorem \ref{main}, we have $|N(s)|\leq \e, \ s \in [-N,N]$. Notice that for
sufficiently small $\e$ we have $h_\e \ll {1\over N}$, and hence the condition $N \leq \frac{1}{2\pi h}$ of Theorem \ref{main} is satisfied. The
measurement results $\Phi(s)$ are identical for $F^0$ and $F^1$, and whatever reconstruction $(\hat A, \hat X)$ of the signal parameters  the algorithm
$\cal R$ produces from $\Phi(s)$, either the distance $d(X^0_l,\hat X_l),$ or $d(X^1_l,\hat X_l),$ is at least $\frac{1}{2} d(X^0_l,X^1_l)$ which, by
Theorem \ref{main}, is not smaller than $\frac{1}{2}C_1 h_\e$. $\square$\par

As for the upper bound on the reconstruction error, we announce the following result:

\bt\label{main1}
Let a signal $F^0=F_{A^0,X^0},$ satisfying assumption $A(m,M)$, form an $(l,h,\rho)$-cluster
$X^0_l=\{x^0_\kappa, x^0_{\kappa+1}, \ldots, x^0_{\kappa+l-1}\}.$ Let the measurements error $\e$ satisfy $\e \leq \e_1=C_3(hN)^{2l}.$ Then solving
the corresponding decimated Prony system of \cite{batenkov_accurate_2014} produces the reconstructed cluster nodes $\bar X^0_l$ with the error at most $\frac{1}{10}\rho h$,
i.e. $d(X^0,\bar X^0)\leq \frac{1}{10}\rho h$.
\et
\par

In particular, since by the assumptions the distance between the cluster nodes is at least $\rho h,$ the number of the nodes, and the inner geometry
of the cluster can be robustly restored.

\par
The proof of Theorem \ref{main1} is based on a combination of the Jacobian estimates in \cite{batenkov_numerical_2014,batenkov_accurate_2014} with the ``Quantitative Inverse Function
theorem'' (Theorem  \ref{Quan.Inverse} below). We plan to present the details separately.

\section{Proof of Theorem \ref{main}}\label{Proof}
\setcounter{equation}{0}

We prove Theorem \ref{main} in several steps. First, for signals $F$ as above we express the Fourier transform ${\cal F}(F)$ through
the moments $m_k(F)$.

\subsection{Fourier transform ${\cal F}(F)$ and moments $m_k(F)$} \label{For.Pro}

For signals $F$ of form (\ref{eq:equation_model_delta}) their Fourier transform ${\cal F}(F)$ can be easily computed explicitly. Let
the moments $m_k(F)$ be defined by

\begin{equation}\label{Moments}
\begin{array}{c}
m_k(F_{A,X})= \int_{-\infty}^{\infty} x^k F_{A,X}(x)dx =\\
= \sum_{j=1}^d a_j x_j^k, \ k=0,1,\ldots .
\end{array}
\end{equation}

\bp\label{expr.F.Tr}
For $F=F_{A,X}=\sum_{j=1}^{d}a_{j}\delta (x-x_{j})$ we have

\be \label{eq:FM.expr}
{\cal F}(F)(s)=\sum_{k=0}^\infty {{m_k(F)}\over {k!}}\tilde s^k, \ \text where \ \ \tilde s = -2\pi i s.
\ee
\ep
\pr

\begin{equation*}
\begin{array}{c}
{\cal F}(F)(s)=\int_{-\infty}^{\infty} e^{-2\pi i sx}F(x)dx= \sum_{j=1}^d a_j e^{-2\pi i x_j s}=\\
= \sum_{j=1}^d a_j \sum_{k=0}^\infty {1\over {k!}}(-2\pi i x_j s)^k = \\
= \sum_{k=0}^\infty {1\over {k!}}(-2\pi i s)^k \sum_{j=1}^d a_jx_j^k = \\
= \sum_{k=0}^\infty {1\over {k!}}m_k(F)\tilde s^k. \ \ \square
\end{array}
\end{equation*}

Thus the Taylor coefficients of the Fourier transform ${\cal F}(F)(s)$ are the consecutive moments $m_k(F)$ divided by $k!$ This fact
provides us an ``Algebraic-Geometric'' approach to the Fourier reconstruction: to produce the signal $F^1$ starting with $F^0$ we analyze
the behavior of the moments $m_k(F)$, and keep them the same for $F^0$ and $F^1$ for $k=0,1,\ldots,2l-2.$ This analysis strongly relies on
recent results in \cite{batenkov_numerical_2014,Bat.Yom1,Bat.Yom2} on the geometry of the ``Prony mapping'', which is formed by the moments $m_k(F)$.

\subsection{Reduction of Theorem \ref{main} to a geometric lemma}

The following result is proved in Section \ref{Geom.Lemma} below:

\bl\label{main.geom}
Let a signal $F^0=F_{A^0,X^0},$ satisfying assumption $A(m,M)$, form an $(d,1,\rho)$-cluster $X^0=\{x^0_1, \ldots, x^0_d\}.$
In other words, all the $d$ nodes of $F^0$ are $\rho$-uniformly distributed in the interval $[-\frac{1}{2},\frac{1}{2}]$. Then there exist
parameters $A^1, X^1$, satisfying assumption $A(\frac{m}{2},2M)$, with $X^1=(x^1_1, \ldots, x^1_d)\subset [-1,1]$, such that\par

\noindent 1. The distance $d(X^0,X^1)$ between the nodes sets $X^0$ and $X^1$ is at least $C_1(m,M,d)>0.$\par

\noindent 2. $m_k(F^0)=m_k(F^1), \ k=0,1,\ldots, 2d-2,$ where $F_1=F_{A^1,X^1}$.
\el
Now we can complete the proof of Theorem \ref{main}. Let a signal $F^0=F_{A^0,X^0},$ satisfying assumption $A(m,M)$, form an $(l,h,\rho)$-cluster
$X^0_l=(x^0_\kappa, x^0_{\kappa+1}, \ldots, x^0_{\kappa+l-1}).$ We rescale the cluster $X^0_l$ to $X'^0_l$ in the interval $[-\frac{1}{2},\frac{1}{2}]$,
with the same amplitudes $A^0_l$. Then we apply Lemma \ref{main.geom} (where we put $d=l$), and find the amplitudes $A^1_l$ and the nodes
$X'^1_l$ in the interval $[-1,1]$, with the same moments as $X'^0_l$ up to $2l-2$. Applying the inverse scaling, and obtain the cluster
$X^1_l \subset [-h,h]$. Clearly, if the moments were equal before shrinking, they will remain equal afterwards. We extend the cluster $X^1_l$ with
the amplitudes $A^1_l$ to the required parameters $A^1, X^1$, adding the non-cluster nodes of $X^0$ with their original amplitudes. By
Lemma \ref{main.geom} the distance $d(X^0_l,X^1_l)$ is at least $C_1h,$ while $m_k(F^0)=m_k(F^1), \ k=0,1,\ldots, 2l-2.$\par

It remains to show that $|{\cal F}(F_{A,X})(s)-{\cal F}(F_{\bar A,\bar X})(s)|\leq C_2 (hs)^{2l-1}.$ By Proposition \ref{expr.F.Tr} we have

\be \label{eq:FM.expr.diff}
{\cal F}(F_{A,X})(s)-{\cal F}(F_{\bar A,\bar X})(s)=\sum_{k=0}^\infty \frac {\gamma_k}{k!}\tilde s^k,
\ee
where $\gamma_k=m_k(F_{A,X})-m_k(F_{\bar A,\bar X}),$ and $\tilde s = -2\pi i s.$ But by our construction $\gamma_k=0, \ k=0,1,\ldots,2l-2.$
On the other hand, since both $X_l$ and $\bar X_l$ are inside $[-h,h]$, while the amplitudes are bounded by $2M,$ we have for any $k$ that
$|\gamma_k|\leq 4lMh^k:=C_3h^k.$ So in fact, for $|s|\leq \frac{1}{2\pi h}$ we get

\begin{equation*}
\begin{array}{c}
\displaystyle
|{\cal F}(F_{A,X})(s)-{\cal F}(F_{\bar A,\bar X})(s)|\leq \sum_{k=2l-1}^\infty \frac {C_3}{k!}|\tilde sh|^k \leq \\
\displaystyle
\leq C_3(2\pi sh)^{2l-1} \sum_{q=0}^\infty \frac {1}{(2l-1+q)!}(2\pi sh)^q \leq \\
\displaystyle
\leq \frac {2C_3(2\pi)^{2l-1}}{(2l-1)!} (sh)^{2l-1}=C_2(sh)^{2l-1},
\end{array}
\end{equation*}
where we put $C_2=\frac {2C_3(2\pi)^{2l-1}}{(2l-1)!}$. This completes the proof of Theorem \ref{main}. $\square$\par

\subsection{Proof of Lemma \ref{main.geom}}\label{Geom.Lemma}

Let a signal $F^0=F_{A^0,X^0},$ satisfying assumption $A(m,M)$, be given, such that all the $d$ nodes of $F$ are $\rho$-uniformly distributed
in the interval $[-\frac{1}{2},\frac{1}{2}]$. We have to show that there exist parameters $A^1,X^1$, satisfying assumption $A(\frac{m}{2},2M)$,
with $X^1_d=(x^1_1, \ldots, x^1_d)\subset [-1,1]$, such that the system of equations

\begin{equation}\label{moments.equ}
\begin{array}{c}
m_k(F^0)=m_k(F^1), \ F^1=F_{A^1,X^1}, \\
k=0,1,\ldots, 2d-2
\end{array}
\end{equation}
is satisfied, while the distance $d(X^0_d,X^1_d)$ is at least $C_1(m,M,d)>0.$ In other words, we have to show that the projection of the set
$Q\subset {\mathbb R}^d \times {\mathbb R}^d$ consisting of $A^1,X^1$ which satisfy (\ref{moments.equ}), is large enough.\par

Let us recall a definition of the Prony mapping $PM:{\mathbb R}^d \times {\mathbb R}^d \to {\mathbb R}^{2d}$, given in \cite{Bat.Yom2}. It
is provided by

\be\label{Prony.Map}
PM(A,X)=\left (m_0(F_{A,X}), \ldots, m_{2d-1}(F_{A,X})\right ) \in {\mathbb R}^{2d}.
\ee
Put $m^0_k=m_k(F_{A,X}).$ Then $PM(A^0,X^0)=(m^0_0,m^0_1,\ldots,m^0_{2d-1}).$ We shall denote by
$$
\begin{array}{c}
\mu=(\mu_0,\ldots,\mu_{2d-1})=\\
=(m_0-m^0_0, \ m_1-m^0_1, \ \ldots, \ m_{2d-1}-m^0_{2d-1})
\end{array}
$$
the coordinates in ${\mathbb R}^{2d}$, with the origin shifted to $\mu^0=PM(A^0,X^0)$. The following result can be easily derived from Theorem 4.5 of
\cite{Bat.Yom1}:

\bt\label{Inv.Prony.Jac}
At each point $(A^0,X^0) \in {\mathbb R}^d \times {\mathbb R}^d$ satisfying conditions of Lemma \ref{main.geom}, the Jacobian $JPM=JPM(A^0,X^0)$ of
the Prony mapping $PM$ is invertible, the norm of its inverse $JPM^{-1}$ is bounded from above by the constant $C_5(m,M,d,\rho)$, and for each
$\mu \in \C^{2d}$ we have $||JPM^{-1}(\mu)||\geq C_6||\mu||,$ with a positive $C_6=C_6(m,M,d,\rho).$\par

Let $L=O\mu_{2d-1}$ be the last coordinate axis in ${\mathbb R}^{2d}$, and let
$$
P_X: {\mathbb R}^d \times {\mathbb R}^d \to {\mathbb R}^d
$$
denote the projection of the signal parameters $A,X$ to the nodes coordinates $X$. Then for each $\mu \in L$ we have
$||P\circ JPM^{-1}(\mu)||\geq C_7||\mu||,$ with a positive $C_7=C_7(m,M,d,\rho).$
\et
\pr The first three statements of Theorem \ref{Inv.Prony.Jac} follow directly from Theorem 4.5 of \cite{Bat.Yom1}. The last statement follows from
the fact that for the fixed nodes the restriction of the first $d$ coordinates of the Prony mapping $PM$, and hence, of its Jacobian $JPM$, to the
coefficients $A$ is a non-degenerate linear mapping to ${\mathbb R}^d$, with the Vandermonde matrix on the nodes $X$. Hence, the pre-image $JPM^{-1}(L)$
cannot be contained in ${\mathbb R}^d \times \{0\} \subset {\mathbb R}^d \times {\mathbb R}^d$. Otherwise at least one of the first $d$ moments would
change along $L$. In fact, it is easy to show that the line $JPM^{-1}(L)$ forms a positive angle with ${\mathbb R}^d \times \{0\}$, which is bounded
from below by a constant depending only on $m,M,d,\rho$. But this is equivalent to the last statement of Theorem \ref{Inv.Prony.Jac}. $\square$\par

On the other hand, using the standard solution procedure of the Prony system, one can easily show the following fact:

\bp\label{Inv.Prony.Jac1}
There are constants $R_1=R_1(m,M,d,\rho)$ and $C_8=C_8(m,M,d,\rho),$ such that for each point $(A^0,X^0) \in {\mathbb R}^d \times {\mathbb R}^d$,
satisfying conditions of Lemma \ref{main.geom}, the first and the second derivatives of the Prony mapping $PM$ are bounded by $C_8$ in the ball
$B_{R_1}$ in ${\mathbb R}^d \times {\mathbb R}^d$ centered at $(A^0,X^0)$.
\ep
With these two preparatory results we now apply the following ``Quantitative Inverse Function Theorem'' (see. e.g. \cite{Yom1}):

\bt\label{Quan.Inverse}
Let $G: (B^m_1,0) \to ({\mathbb R}^m,0)$ be a twice differentiable mapping of the unit ball $B^m_1$ at the origin in ${\mathbb R}^m$ to
${\mathbb R}^m$, with $G(0)=0$, such that the Jacobian $J=JG(0)$ is invertible, and $||J^{-1}||\leq K_1 <\infty.$ Assume that the second derivatives
of $G$ are bounded by $K_2$ in the ball $B^m_1.$ Then the inverse mapping $G^{-1}$ exists in the ball $B^m_{R_2}$ of radius $R_2$, centered at
$0\in {\mathbb R}^m$, and satisfies there the condition
\be\label{eq:Inv.Prony.Jac}
\begin{array}{c}
||G^{-1}(x)- J^{-1}(x)||\leq C_9(m,K_1,K_2) ||x||^2, \\
\text for  \ \ x\in B^m_{R_2},
\end{array}
\ee
with $R_2=R_2(m,K_1,K_2)$ and $C_9=C_9(m,K_1,K_2)$ depending only on $m, K_1,$ and $K_2$.
\et
We apply Theorem \ref{Quan.Inverse}, properly rescaled, to the the Prony mapping $PM$ in the ball $B_{R_1}$ in ${\mathbb R}^d \times {\mathbb R}^d$
centered at $(A^0,X^0)$, with the bounds provided by Theorem \ref{Inv.Prony.Jac} and Proposition \ref{Inv.Prony.Jac1}. We conclude that in the ball
$B_{R_3}$ of radius $R_3$ at the point $\mu^0=PM(A^0,X^0)\in {\mathbb R}^{2d}$ the inverse $PM^{-1}$ exists and satisfies

\be\label{eq:Inv.Prony.Jac1}
||PM^{-1}(\mu)- J^{-1}(\mu)||\leq C_{10} ||\mu||^2,
\ee
where $J=JPM$ is the Jacobian of the Prony mapping $PM$ at $(A^0,X^0)$, and the constants $R_3$ and $C_{10}$ depend only on $m,M,d,\rho$.\par

From the last statement of Theorem \ref{Inv.Prony.Jac} we get for $\mu \in L=O\mu_{2d-1}$

\be\label{lower.b}
||P\circ J^{-1}(\mu)||\geq C_{11}||\mu||,
\ee
with $P$ the projection of the signal parameters $A,X$ to the nodes $X$, and $C_{11}=C_{11}(m,M,d,\rho)$ a positive constant.\par

Finally we put $\mu^1=(0,\ldots,0,\eta)$, with $\eta = \min\bigl(R_3,\frac {C_{11}}{2C_{10}}\bigr)$, and take $(A^1,X^1)$ to be the inverse image $PM^{-1}(\mu^1)$. By
the construction we have

\be\label{moments.equ1}
m_k(F^0)=m_k(F^1), \ k=0,1,\ldots, 2d-2.
\ee
On the other hand, $X^1=P(PM^{-1}(\mu^1))=P(J^{-1}(\mu^1)-w),$ where $w=J^{-1}(\mu^1)-PM^{-1}(\mu^1)$, and hence, by (\ref{eq:Inv.Prony.Jac1}),
we have $||w||\leq C_{10} ||\mu^1||^2=C_{10}\eta^2.$ By (\ref{lower.b}) we get $||P\circ J^{-1}(\mu^1)||\geq C_{11}\eta$, and therefore

$$
\begin{array}{c}
||X^1||=||P(J^{-1}(\mu^1)-w)||\geq C_{11}\eta - C_{10}\eta^2=\\
=\eta(C_{11} - C_{10}\eta) \geq\eta\frac {C_{11}}{2}:=C_{12}.
\end{array}
$$

The norm $||X^1||$ of $X^1$ here is the $l^2$ norm with respect to the coordinates in ${\mathbb R}^d$ centered at $X^0$. Hence $\frac{1}{\sqrt l}||X^1||$
bounds from below the distance $d(X^0,X^1)$. This completes the proof of Lemma \ref{main.geom}, with $C_1= \frac{1}{\sqrt l}C_{12}$. $\square$.

\smallskip

\noindent \textbf{Remark.} In this paper we consider only the curve $PM^{-1}(L)$ where the first $2l-2$ moments $m_k$ take equal value. In the direction of this
curve the magnification of the
measurements error is maximal. In fact, for each $q=1,\ldots,2l-2$ there is a stratum $\Sigma_q$ in ${\mathbb R}^d \times {\mathbb R}^d$, of dimension
$2l-q-1$, where the first $q$ moments $m_k$ take equal value. In the direction of this stratum the error magnification is of order $q+1$. The geometry of
the strata $\Sigma_q$ plays important role in the understanding of the error magnification patterns which occur in the Fourier reconstruction of
spike-trains. We plan to present the results in this direction separately.

\section*{Some examples}

The following examples illustrate the shape and behavior of the signals $F^0_q$ and $F^1_q, \ q=1,3,5,$ for which the difference
$DF_q(s)={\cal F}(F^0_q)(s)-{\cal F}(F^1_q)(s)$ between their Fourier transforms is of order $q$ in $s$. As it was explained above, the
geometry of the strata $\Sigma_q$, containing $F_q$ plays important role in the error magnification which occurs in the Fourier reconstruction
of spike-trains.

\smallskip

We consider signals with $d=3$ nodes of the form \eqref{eq:equation_model_delta}: $F_q(x)=\sum_{j=1}^3 a_{qj}\delta (x-x_{qj})$. Their specific
parameters are shown in table \ref{tbl:experiment-signals}. In this table we assume $h$ to be fixed, and put $\eta=\tilde \eta h$, with $\tilde \eta$
being the ``free parameter along the stratum $\Sigma_q$''. The maximal distance between the nodes of $F^0_q$ and $F^1_q$ in each case is $2\eta$.

\smallskip

\begin{figure}[!h]
\centering
\includegraphics[width=3in]{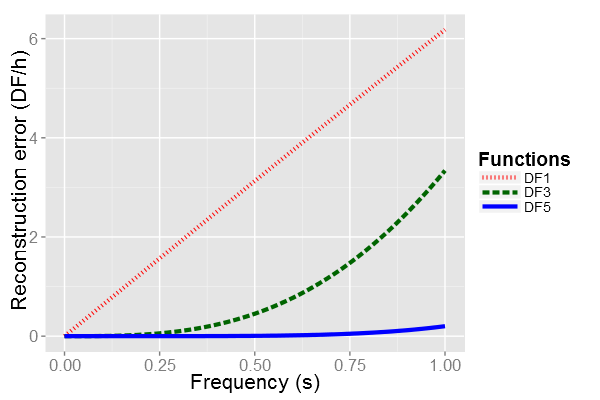}
\caption{Fourier differences $DF$ of the signals from Table \ref{tbl:experiment-signals}.}
\label{img:example-plot1}
\end{figure}

Table \ref{tbl:experiment-moments} shows the difference between the moments $m_0, m_1, m_2, m_3,m_4$ of $F^0_q$ and $F^1_q$. We see that this
difference is zero exactly for the first $q$ moments, $q=1,3,5$. Figure \ref{img:example-plot1} shows the difference
$DF_q(s)={\cal F}(F^0_q)(s)-{\cal F}(F^1_q)(s)$, for the frequency $s \in [0,1]$. The normalized difference presented on the ordinate is $\frac{DF}{h}$.
In this figure we fix $\eta=0.05, h=0.1 $.

\begin{table}[!h]
\caption{Signals parameters}
\label{tbl:experiment-signals}
\centering
\begin{tabular}{|c||c|c|c||c|c|c|}
\hline
~ & $a_1$ & $a_2$ & $a_3$ & $x_1$ & $x_2$ & $x_3$ \\
\hline
\hline
$F^0_1$ & 1 & 1 & 1 & $-h-\eta$ & $-\eta$ & $+h+\eta$ \\
\hline
$F^1_1$ & 1 & 1 & 1 & $-h-\eta$ & $+\eta$ & $+h+\eta$ \\
\hline
\hline
$F^0_3$ & 1 & 1 & 1 & $-h-\eta$ & $-\eta$ & $+h+2\eta$ \\
\hline
$F^1_3$ & 1 & 1 & 1 & $-h-2\eta$ & $+\eta$ & $+h+\eta$ \\
\hline
\hline
$F^0_5$ & $-1-3\tilde\eta$ & $2+3\tilde\eta$ & $-1$ & $-h-\eta$ & $-\eta$ & $h+2\eta$ \\
\hline
$F^1_5$ & $-1$ & $2+3\tilde\eta$ & $-1-3\tilde\eta$ & $-h-2\eta$ & $+\eta$ & $+h+\eta$ \\
\hline
\end{tabular}
\end{table}

\begin{table}[!h]
\caption{Moments differences}
\label{tbl:experiment-moments}
\centering
\begin{tabular}{|c||c|c|c|c|c|}
\hline
~ & $\Delta m_0$ & $\Delta m_1$ & $\Delta m_2$ & $\Delta m_3$ & $\Delta m_4$ \\
\hline
\hline
$F_1$  & 0 & $2\eta$ & 0 & $2\eta^3$ & 0\\
\hline
\hline
$F_3$  & 0 & 0 & 0 & $6h^2\eta+18h\eta^2+16\eta^3$ & 0\\
\hline
\hline
$F_5$  & 0 & 0 & 0 & 0 & 0\\
\hline
\end{tabular}
\end{table}

\section*{Acknowledgment}
This research was supported by the ISF Grant No. 779/13.


\end{document}